\documentclass{gtart_h}  

\def\ifplaintex{\expandafter\ifx\csname documentclass\endcsname\relax}

\ifplaintex 
\else
\fi

\expandafter\ifx\csname beginpicture\endcsname\relax
\expandafter\ifx\csname documentclass\endcsname\relax
\input pictex \else
\input prepictex \input pictex \input postpictex \fi\fi

\def\gt{{\mathsurround=0pt\it $\cal G\mskip-2mu$eometry \&\ 
$\cal T\!\!$opology}}        %  journal title in recommended style

\def\gtp{{\mathsurround=0pt\it $\cal G\mskip-2mu$eometry \&\ 
$\cal T\!\!$opology $\cal P\!$ublications}}  % GT publications

\def\lognumber#1{\def\thelognumber{#1}}
\def\volumenumber#1{\def\thevolumenumber{#1}}
\def\papernumber#1{\def\thepapernumber{#1}}
\def\volumeyear#1{\def\thevolumeyear{#1}}

\def\pagenumbers#1#2{\def\startpage{#1}\def\finishpage{#2}}
\def\published#1{\def\publishdate{#1}}
\def\proposed#1{\def\theproposer{#1}}
\def\seconded#1{\def\theseconders{#1}}
\def\received#1{\def\receiveddate{#1}}

\def\accepted#1{\def\accepteddate{#1}}
\def\asciititle#1{\def\theasciititle{#1}}
\def\covertitle#1{\def\thecovertitle{#1}}
\def\coverauthors#1{\def\thecoverauthors{#1}}
\def\asciiauthors#1{\def\theasciiauthors{#1}}
\def\asciiaddress#1{\def\theasciiaddress{#1}}
\def\asciiemail#1{\def\theasciiemail{#1}}

\long\def\asciiabstract#1{\long\def\theasciiabstract{#1}}
\def\asciikeywords#1{\def\theasciikeywords{#1}}

\let\\\par\let\thelognumber\relax
\let\thevolumenumber\relax\let\thepapernumber\relax
\let\thevolumeyear\relax\let\thesamplenumber\relax\let\startpage\relax
\let\finishpage\relax\let\publishdate\relax\let\receiveddate\relax
\let\reviseddate\relax\let\accepteddate\relax\let\theasciititle\relax
\let\thecovertitle\relax\let\theasciiauthors\relax\let\theasciiaddress\relax
\let\theasciiabstract\relax\let\theasciikeywords\relax
\let\theasciiemail\relax\let\theshortauthors\relax\let\theshorttitle\relax
\let\thecoverauthors\relax

\long\def\maketitlep{   % start of definition of \maketitlep

\count0=\startpage

\gt\hfill      %   Journal title (top left) 
\beginpicture
\setcoordinatesystem units <0.33truein, 0.33truein> point at 2.2 0.9
\setplotsymbol ({$\cal G$})
\plotsymbolspacing=9truept
\circulararc 315 degrees from 0 1 center at 0 0
\setplotsymbol ({$\cal T$})
\circulararc 315 degrees from 1 -1 center at 1 0
\endpicture
\break
{\small\ifx\thesamplenumber\relax % sample?  
Volume \else Sample
\fi\thevolumenumber\ (\thevolumeyear)
\startpage--\finishpage\nl
Published: \publishdate}
\vglue 0.5truein plus 0.4fil minus 0.1truein

{\parskip=0pt\leftskip 0pt plus 1fil\def\\{\par\smallskip}{\ifplaintex\large
\else\Large\fi\bf\thetitle}\par\medskip}   

\vglue 0pt plus 0.1fil 

{\parskip=0pt\leftskip 0pt plus 1fil\def\\{\par}{\sc\theauthors}
\par\medskip}

\vglue 0pt plus 0.1fil 

{\small\parskip=0pt\let\newline\\
{\leftskip 0pt plus 1fil\def\\{\par}{\sl\theaddress}\par}
\expandafter\ifx\theemail\relax    % email address?
\relax\else\vglue 5pt plus 0.02fil minus 2pt\def\\{\stdspace{\rm 
and}\stdspace} 
\cl{Email:\stdspace\tt\theemail}\fi
\ifx\theurl\relax                  % URL given?
\relax\else\vglue 5pt plus 0.02fil minus 2pt\def\\{\stdspace{\rm 
and}\stdspace}
\cl{URL:\stdspace\tt\theurl}\fi\par}

\vglue 7pt plus 0.3fil minus 3pt

{\bf Abstract}
\vglue 5pt plus 0.1fil minus 2pt

\theabstract

\vglue 7pt plus 0.3fil minus 3pt

{\bf AMS Classification numbers}\quad Primary:\quad \theprimaryclass

Secondary:\quad \thesecondaryclass

\vglue 5pt plus 0.3fil minus 2pt

{\bf Keywords:}\quad \thekeywords

\vglue 10pt plus 0.5fil minus 5pt

{\small  Proposed: \theproposer\hfill Received: \receiveddate\nl
Seconded: \theseconders\hfill 
\ifx\reviseddate\relax                         % paper revised?
Accepted: \accepteddate                        % no
\else
Revised: \reviseddate                          % yes
\fi}
\eject
}       %  end of definition of \maketitlep

\let\maketitlepage\maketitlep
\let\maketitle\maketitlepage

\font\phead=cmsl9 scaled 950
\font=cmsl9 scaled 1050
\font\pnum=cmbx10 scaled 913
\font\lnum=cmbx10 
\font\pfoot=cmsl9 scaled 950
\font=cmsl9 scaled 1050
\ifplaintex
\headline{\vbox to 0pt{\vskip -4.5mm\line{\small\phead\ifnum
\count0=\startpage ISSN 1364-0380 (on line)
1465-3060 (printed) \hfill {\pnum\folio}\else\ifodd\count0\def\\{ }% 
\ifx\theshorttitle\relax\thetitle\else\theshorttitle\fi\hfill{\pnum\folio}
\else\def\\{ and }{\pnum\folio}\hfill\ifx\theshortauthors\relax\theauthors
\else\theshortauthors\fi\fi\fi}\vss}}
\footline{\vbox to 0pt{\vglue 0mm\line{\small\pfoot\ifnum\count0=\startpage
\copyright\ \gtp\hfill\else
\gt, Volume \thevolumenumber\ (\thevolumeyear)\hfill\fi}\vss
}}
\else
\makeatletter
\def\@oddhead{{\small\lhead\ifnum\count0=\startpage ISSN 1364-0380 (on line)
1465-3060 (printed) \hfill {\lnum\number\count0}\else\ifodd\count0
\def\\{ }\ifx\theshorttitle\relax \thetitle \else\theshorttitle\fi\hfill
{\lnum\number\count0}\else\def\\{ and }{\lnum\number\count0}
\hfill\ifx\theshortauthors\relax 
\theauthors\else\theshortauthors\fi\fi\fi}}\def\@evenhead{\@oddhead}
\def\@oddfoot{\small\lfoot\ifnum\count0=\startpage\copyright\ \gtp\hfill\else
\gt, Volume \thevolumenumber\ (\thevolumeyear)\hfill\fi}
\def\@evenfoot{\@oddfoot}
\makeatother
\fi

\newwrite\gtoutfile
\long\gdef\makeheadfile{  %%% start of definition of \makeheadfile
{\def\\{, }\def\s{ }
\immediate\openout\gtoutfile head.xxx
\immediate\write\gtoutfile{Proxy-for: \ifx\theasciiauthors\relax
\theauthors\else\theasciiauthors\fi\s<\ifx\theasciiemail\relax\theemail\else\theasciiemail\fi>}
\immediate\write\gtoutfile{\noexpand\\}
\immediate\write\gtoutfile{Authors: \ifx\theasciiauthors\relax
\theauthors\else\theasciiauthors\fi}
{\def\\{ }\immediate\write\gtoutfile{Title: \ifx\theasciititle\relax
\thetitle\else\theasciititle\fi}}
\immediate\write\gtoutfile{Subj-class: GT or SG or MG etc}
\immediate\write\gtoutfile{MSC-class: \theprimaryclass\ifx\thesecondaryclass\relax\else, \thesecondaryclass\fi}
\immediate\write\gtoutfile{Journal-ref: Geom. Topol. \thevolumenumber
(\thevolumeyear) \startpage-\finishpage}
\immediate\write\gtoutfile{Comments: Published by Geometry and Topology at}
\immediate\write\gtoutfile{\s\s http://www.maths.warwick.ac.uk/gt/GTVol\thevolumenumber/paper\thepapernumber.abs.html}
\immediate\write\gtoutfile{\noexpand\\}
\immediate\write\gtoutfile{}
\ifx\theasciiabstract\relax
\immediate\write\gtoutfile{\theabstract}\else
\immediate\write\gtoutfile{\theasciiabstract}\fi
\immediate\write\gtoutfile{}
\immediate\write\gtoutfile{\noexpand\\}
\immediate\write\gtoutfile{}
\immediate\closeout\gtoutfile}}  %%% end of definition of \makeheadfile

\def\maketitlepage{\maketitlep\makeheadfile}
\let\maketitle\maketitlepage

\lognumber{536}
\received{6 December 2004}
\volumenumber{9}\papernumber{19}\volumeyear{2005}
\pagenumbers{813}{832}   
\published{18 May 2005}
\accepted{2 May 2005}
\proposed{Peter Ozsvath}
\seconded{Ronald Fintushel, Tomasz Mrowka}

\usepackage{amsmath, amssymb, graphicx}

\newtheorem{thm}{Theorem}[section]

\newtheorem{cor}[thm]{Corollary}
\newtheorem{lem}[thm]{Lemma}
\newtheorem{prop}[thm]{Proposition}

\theoremstyle{definition}
\newtheorem{defn}[thm]{Definition}

\newtheorem{rem}[thm]{Remark}

\numberwithin{equation}{section}

\newcommand{\bfz}{{\mathbb {Z}}}
\newcommand{\bfq}{{\mathbb {Q}}}

\newcommand{\s}{\mathbf s}

\newcommand{\C}{\mathbb C}
\newcommand{\Z}{\mathbb Z}

\newcommand{\bfr}{\mathbb R}

\newcommand{\cpkk}{{\overline {{\mathbb C}{\mathbb P}^2}}}
\newcommand{\cpk}{{\mathbb {CP}}^2}
\newcommand{\cphat}{{\mathbb {CP}}^2\# 6{\overline {{\mathbb C}{\mathbb P}^2}}}
\newcommand{\cphet}{{\mathbb {CP}}^2\# 7{\overline {{\mathbb C}{\mathbb P}^2}}}
\newcommand{\cpnyolc}{{\mathbb {CP}}^2\# 8{\overline {{\mathbb C}{\mathbb P}^2}}}
\newcommand{\eegy}{{\mathbb {CP}}^2\# 9{\overline {{\mathbb C}{\mathbb P}^2}}}

\DeclareMathOperator{\III}{III}

\begin{document}

\title{An exotic smooth structure on $\cphat$}
\asciititle{An exotic smooth structure on CP^2+6CP^2-bar}
\covertitle{An exotic smooth structure on ${\noexpand\bf CP}^2\noexpand\#6\noexpand\overline{{\noexpand\bf CP}^2}$}

\author{Andr\'{a}s I Stipsicz\\Zolt\'an Szab\'o}
\asciiauthors{Andras I Stipsicz and Zoltan Szabo}
\coverauthors{Andr\noexpand\'{a}s I Stipsicz\\Zolt\noexpand\'an Szab\noexpand\'o}

\address{R\'enyi Institute of Mathematics, Hungarian Academy of Sciences\\
H-1053 Budapest, Re\'altanoda utca 13--15, Hungary\\
{\rm and}\\
Institute for Advanced Study, Princeton, NJ 08540, USA\\\smallskip\\
{\rm Email:\qua}{\tt \mailto{stipsicz@renyi.hu}, \mailto{stipsicz@math.ias.edu}}
\\\medskip\\Department of Mathematics, Princeton University\\
Princeton, NJ 08544, USA\\\smallskip\\
{\rm Email:\qua}{\tt \mailto{szabo@math.princeton.edu}}}

\asciiaddress{Renyi Institute of Mathematics, Hungarian Academy of Sciences\\
H-1053 Budapest, Realtanoda utca 13--15, Hungary\\
Institute for Advanced Study, Princeton, NJ 08540, 
USA\\and\\Department of Mathematics, Princeton University\\
Princeton, NJ 08544, USA}

\asciiemail{stipsicz@renyi.hu, stipsicz@math.ias.edu and
szabo@math.princeton.edu}

\begin{abstract}
We construct smooth 4--manifolds homeomorphic but not
diffeomorphic to $\cphat$.
\end{abstract}
\asciiabstract{We construct smooth 4-manifolds homeomorphic but not
diffeomorphic to CP^2+6CP^2-bar.}
\asciikeywords{Exotic smooth 4-manifolds, Seiberg-Witten invariants,
rational blow-down, rational surfaces}

\maketitle

\section{Introduction}
Based on work of Freedman \cite{Fr} and Donaldson \cite{D}, in the mid
80's it became possible to show the existence of exotic smooth
structures on closed simply connected 4--manifolds.  On one hand,
Freedman's classification theorem of simply connected, closed
topological 4--manifolds could be used to show that various
constructions provide homeomorphic 4--manifolds, while the computation
of Donaldson's instanton invariants provided a smooth invariant
distinguishing appropriate examples up to diffeomorphism, see
\cite{D1} for the first such computation. For a long time the pair
$\cpnyolc$ (the complex projective plane blown up at eight points) and
a certain algebraic surface (the Barlow surface) provided such a
simply connected pair with smallest Euler characteristic \cite{Kot}.
Recently, by a clever application of the rational blow-down operation
originally introduced by Fintushel and Stern \cite{FS1}, Park found a
smooth 4--manifold homeomorphic but not diffeomorphic to $\cphet$
\cite{P}. Applying a similar rational blow-down construction we show
the following:
\begin{thm}\label{t:main}
There exists a smooth 4--manifold $X$ which is homeomorphic to $\cphat$
but not diffeomorphic to it.
\end{thm}

Note that $X$ has Euler characteristic $\chi (X)=9$, and thus
provides the smallest known closed exotic simply connected smooth 4--manifold.
The proof of Theorem~\ref{t:main} involves two steps. First we will
construct a smooth 4--manifold $X$ and  determine its fundamental
group and characteristic numbers. Applying Freedman's theorem, we conclude
that $X$ is homeomorphic to $\cphat$. Then by computing the Seiberg--Witten
invariants of $X$ we show that it is not diffeomorphic to $\cphat$.
By determining all Seiberg--Witten basic classes of $X$ we can also show that
it is minimal. This result, in conjunction with the result
of \cite{OSZP} gives:

\begin{cor}\label{hanykul}
Let $n \in \{ 6,7,8\}$. Then there are at least $n-4$ different smooth
structures on the topological manifolds $\cpk \# n \cpkk$.  The
different smooth 4--manifolds $Z_1(n), Z_2(n), \ldots ,Z_{n-4}(n)$
homeomorphic to $\cpk \# n \cpkk$ have $0,2,\ldots , 2^{n-5}$
Seiberg--Witten basic classes, respectively. \qed
\end{cor}

In Section~\ref{s:top} we give several  constructions of
exotic smooth structures on the topological 4--manifold
$\cphat$ by rationally blowing down
various configurations of chains of 2--spheres. Since the generalized
rational blow-down operation is symplectic when applied along
symplectically embedded spheres (see \cite{Sym}), the 4--manifolds that
are constructed here  all admit symplectic structures.  
The computation of their Seiberg--Witten basic classes show that they 
are all minimal symplectic 4--manifolds with isomorphic Seiberg--Witten
invariants. It is not known whether these examples are diffeomorphic
to each other. 

It is interesting to note that any two  minimal symplectic 4--manifolds 
 on the topological manifold  $\cpk
\# n \cpkk$  $n\in \{ 1, \ldots , 8\}$ have (up to sign)
identical Seiberg--Witten invariants. As a 
corollary, Seiberg--Witten invariants can tell apart only at most finitely
many symplectic structures on the topological manifold  $\cpk
\# n \cpkk$ with $n\leq 8$.

{\bf Acknowledgements}\qua We would like to thank Andr\'as N\'emethi,
Peter Ozsv\'ath, Jongil Park and Ron Stern for enlightening
discussions.  The first author was partially supported by OTKA T49449
and the second author was supported by NSF grant number DMS 0406155.%
\footnote{After the submission of this paper the results of 
Theorem~\ref{t:main} and Corollary~\ref{hanykul} have been improved
by finding infinitely many exotic smooth structures on 
$\cpk \# n \cpkk$ for $n\geq 5$, see \cite{FSuj, PSS}.}

\section{The topological constructions}
\label{s:top}

In constructing the 4--manifolds encountered in Theorem~\ref{t:main}
we will apply the generalized rational blow-down operation
\cite{Pratb} to certain configurations of spheres in rational
surfaces. In order to locate the particular configurations, we start
with a special elliptic fibration on $\eegy$. The proof of the
following proposition is postponed to Section~\ref{app}.  (For conventions
and constructions see \cite{HKK}.)

\begin{figure}[htb]
\begin{center}
\setlength{\unitlength}{1mm}
\includegraphics[height=2cm]{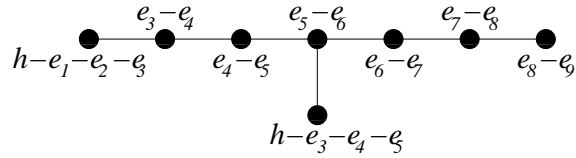}
\end{center}
\caption{Plumbing diagram of the singular fiber of type $\III^*$}
\label{f:III}
\end{figure}

\begin{prop}\label{hkk}
There is an elliptic fibration $\eegy \to {\mathbb {CP}}^1$ with a singular
fiber of type $\III^*$, three fishtail fibers and two sections.
\qed
\end{prop}
The type $\III^*$ singular fiber (also known as the ${\tilde {E}}_7$
singular fiber) can be given by the plumbing diagram of
Figure~\ref{f:III}.  (All spheres in the plumbing have
self--intersection equal to $-2$.)
If $h, e_1, \ldots , e_9$ is the standard generating system of
$H_2(\eegy ; \Z )$ then the elliptic fibration can be 
arranged so that the homology classes of the spheres in the
$\III^*$ fiber are equal to the classes given in
Figure~\ref{f:III}. We also show in Section~\ref{app} that the two 
sections can be chosen to intersect the spheres in the left and the right ends
of Figure~\ref{f:III}, respectively.

\subsection{Generalized rational blow-down}
Let $L_{p,q}$ denote the lens space $L(p^2, pq-1)$, where $p\geq q\geq
1$ and $p,q$ are relatively prime. Let $C_{p,q}$ denote the plumbing 4--manifold
obtained by plumbing 2--spheres along the linear graph with decorations
$d_i\leq -2$ given by the continued fractions of $-\frac{p^2}{pq-1}$;
we have the obvious relation $\partial C_{p,q}=L_{p,q}$, cf also \cite{Pratb}.  Let
$K\in H^2 (C_{p,q}; \Z )$ denote the cohomology class which evaluates
on each 2--sphere of the plumbing diagram as $d_i +2$.

\begin{prop}{\rm\cite{CH, Pratb, Sym}}\label{p:ch}\qua 
The 3--manifold $\partial C_{p,q}=L(p^2, pq-1)$ bounds a rational 
ball $B_{p,q}$ and the cohomology class $K\vert _{\partial C_{p,q}}$ extends 
to $B_{p,q}$. \qed
\end{prop}

The following proposition provides embeddings of some of the above
plumbings into rational surfaces.

\begin{prop}
\begin{itemize}
\item
The 4--manifold $C_{28,9}$ embeds into $\cpk \# 17 \cpkk$; 
\item 
$C_{46,9}$ embeds into $\cpk \# 19 \cpkk$, and finally
\item 
$C_{64,9}$ embeds into $\cpk \# 21 \cpkk$.
\end{itemize}
\end{prop}
\begin{rem}
The linear plumbings giving the configurations considered above 
are as follows:
\begin{itemize}
\item $C_{28,9}=(-2,-2,-12,-2,-2,-2,-2,-2,-2,-2,-4)$,
\item $C_{46,9}=(-2,-2,-2,-2,-12,-2,-2,-2,-2,-2,-2,-2,-6)$ and
\item $C_{64,9}=(-2,-2,-2,-2,-2,-2,-12,-2,-2,-2,-2,-2,-2,-2,-8)$.
\end{itemize}
\end{rem}
\begin{proof}
Let us consider an elliptic fibration on $\eegy$ with a type $\III^*$
singular fiber, three fishtail fibers $F_1, F_2, F_3$ and two sections
$s_1, s_2$ as described by the schematic diagram of
Figure~\ref{f:singfib}.  Let $A_i$ denote the intersection of $F_i$
with the section $s_2$, while $B_i$ denotes the intersection of the
fiber $F_i$ with $s_1$ ($i=1,2,3)$.
\begin{figure}[htb]
\begin{center}
\setlength{\unitlength}{1mm}
\includegraphics[height=5.5cm]{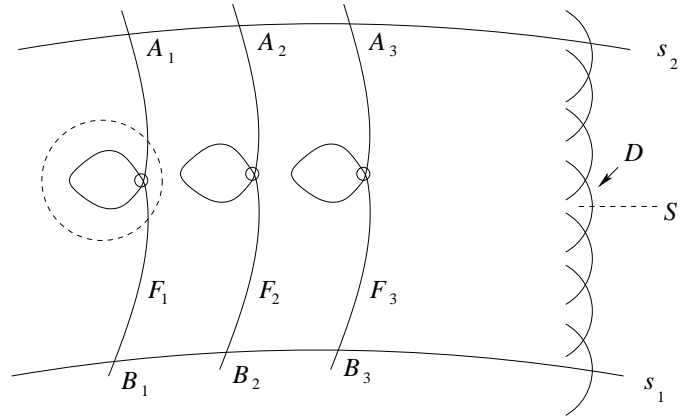}
\end{center}
\caption{Singular fibers in the fibration}
\label{f:singfib}
\end{figure}
First blow up the three double points 
(indicated by small circles) of the three fishtail fibers.
To get the first configuration, further blow up at $A_1, A_2, A_3$ 
and smooth the transverse intersections $B_1, B_2, B_3$. Finally, apply two
more blow-ups inside the dashed circle as shown by Figure~\ref{f:inside}. 
\begin{figure}[htb]
\begin{center}
\setlength{\unitlength}{1mm}
\includegraphics[height=7cm]{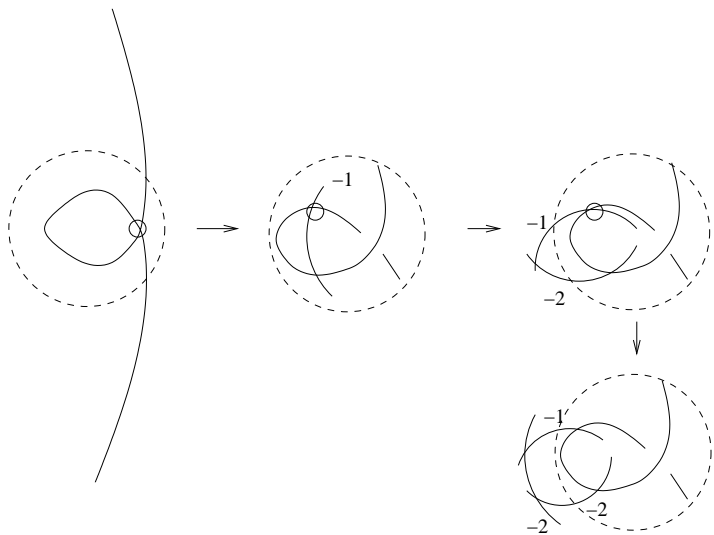}
\end{center}
\caption{Further blow-ups of the fishtail fiber}
\label{f:inside}
\end{figure}
By counting the number of blow-ups, the desired embedding of
$C_{28,9}$ follows.

In a similar way, now blow up $A_1, A_2, B_3$, and smooth
$B_1, B_2$ and $A_3$. Four further blow-ups in the manner 
depicted by Figure~\ref{f:inside} provides the embedding 
of $C_{46,9}$. 

Finally, by blowing up $A_1, B_2, B_3$, and smoothing $B_1, A_2$ and
$A_3$, and then performing six further blow-ups as before inside the
dashed circle, we get the embedding of $C_{64,9}$ as claimed.
\end{proof}

\begin{lem}
For $i=0,1,2$ the embedding
$C_{28+18i,9}\subset \cpk \# (17+2i)\cpkk$ found above
has simply connected complement.
\end{lem}
\begin{proof}
Since rational surfaces are simply connected, the simple connectivity
of the complement follows once we show that a circle in the boundary
of the complement is homotopically trivial. Recall that, since the
boundary of the complement is a lens space, it has cyclic fundamental
group. In conclusion, homotopic triviality needs to be checked only
for the generator of the fundamental group of the boundary.  We claim
that the normal circle to the $(-2)$--framed sphere $D$ in the $\III^*$
fiber intersected by the dashed $(-2)$--curve $S$ of
Figure~\ref{f:singfib} (which is in the $\III^*$ fiber but not in our
chosen configuration) is a generator of the fundamental group of the
boundary 3--manifold. This observation easily follows from the
facts that for the boundary lens space the first homology is naturally
isomorphic to the fundamental group, and in the first homology the
normal circle in question is $13+8i$ times the generator given by the
linking normal circle of the last sphere of the configuration
($i=0,1,2$). Since for $i=0,1,2$ we have that $13+8i$ is relatively
prime to $28+18i$, and the hemisphere of $S$ in the complement shows
that the normal circle of $D$ is homotopically trivial in the
complement, the proof of the lemma follows.
\end{proof}

\begin{rem}
It is not hard to show that the 3--manifold $\partial C_{28,9}$ does
bound a rational ball: we can embed $C_{28,9}$ into $11\cpkk$ and the
closure of the complement of the embedding (with reversed orientation)
can be easily seen to be an appropriate rational ball. In turn,
the embedding $C_{28,9}\subset 11\cpkk$ results from the
following observation. Attach a 4--dimensional 2--handle to $C_{28,9}$
along the $(-1)$--framed unknot $K$ indicated by the plumbing diagram
of Figure~\ref{f:plum}.
\begin{figure}[htb]
\begin{center}
\setlength{\unitlength}{1mm}
\includegraphics[height=3cm]{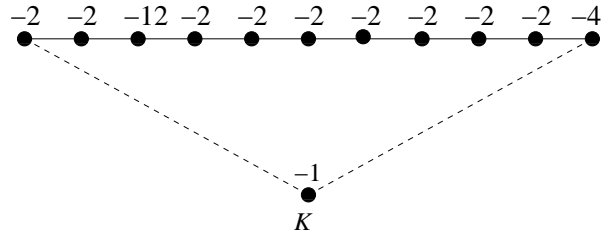}
\end{center}
\caption{Plumbing diagram of the 4--manifold $C_{28,9}$} 
\label{f:plum}
\end{figure}
By subsequently sliding down the $(-1)$--framed unknots we arrive at a
0--framed unknot, showing that the handle attachment along $K$ embeds
$C_{28,9}$ into a 4--manifold diffeomorphic to the connected sum
$S^2\times D^2 \# 11 \cpkk$.  By attaching a 3-- and a 4--handle to
this 4--manifold we get a closed 4--manifold diffeomorphic to
$11\cpkk$. The appropriate modification of the procedure gives the
rational balls $B_{46,9}$ and $B_{64,9}$.  In fact, by adding a
cancelling 1--handle to $K$, doing surgery along it, following the
resulting $0$--framed unknot during the blow-downs of the
$(-1)$--spheres as described above, and then doing another surgery
along the resulting 0--framed circle, we arrive at an explicit surgery
description of the 4--manifold $B_{28,9}$ (and similarly of $B_{46,9}$
and $B_{64,9}$).  Notice that this procedure also shows that the
rational homology balls $B_{28+18i,9}$ we get in this way admit
handlebody decompositions involving only handles in dimensions 0,1
and 2, hence the maps $\pi _1(\partial B_{28+18i,9})\to \pi
_1(B_{28+18i,9})$ induced by the natural embeddings are surjective.
We just note here that the same argument works for all linear
plumbings $(b_1, \ldots , b_k)$ with $b_i\leq -2$ ($i=1,\ldots , k)$
we get from the plumbing $(-4)$ by the repeated applications of the
following two transformation rules (cf also \cite{KS}):
\begin{itemize}
\item $(b_1, \ldots , b_k )\longrightarrow (b_1-1, b_2, \ldots , b_k, -2)$ and
\item $(b_1, \ldots , b_k )\longrightarrow (-2, b_1, \ldots ,b_{k-1}, b_k-1)$.
\end{itemize}
\end{rem}

We will give the details of the computation of Seiberg--Witten
invariants in Section~\ref{s:sw} only for the rational blow-down of
$C_{28,9}\subset \cpk \# 17 \cpkk$.  To make this computation explicit,
we fix the convention that the second homology group $H_2(\cpk \# 17
\cpkk ; \Z)$ is generated by the homology elements $h, e_1, \ldots ,
e_{17}$ (with $h^2=1, e_i^2=-1$ $i=1,\ldots , 17$) and in this basis
the homology classes of the spheres in $C_{28,9}$ can be given (from
left to right on the linear plumbing of Figure~\ref{f:plum}) as
\begin{gather*}
e_{16}-e_{17}, \ \ e_{10}-e_{16},\\ 
9h-2e_1- \sum _3 ^9 3e_i-\sum _{10}^{12} 2e_i- \sum _{13} ^{17} e_i, \ \ \ \
 h-e_1-e_2-e_3\\
 e_3-e_4, \  e_4-e_5, \ e_5-e_6, \  e_6-e_7, \ e_7-e_8, \ e_8-e_9, \ 
e_9-e_{13}-e_{14}-e_{15},
\end{gather*}
where here the two sections $s_1$, $s_2$ represent $e_1$ and $e_9$,
respectively.  

\begin{defn}\label{d:ex}
Let us define $X_1$ as the rational blow-down of $\cpk \# 17 \cpkk$
along the copy of $C_{28,9}$ specified above, that is,
\[
X_1=(\cpk \# 17 \cpkk - {\mbox { int }}(C_{28,9}))\cup _{L(784,251)}B_{28,9} .
\]
Similarly, $X_2, X_3$ is defined as the rational blow-down of the configurations
$C_{46,9}$ and $C_{64,9}$ in the appropriate rational surfaces.
\end{defn}
As a consequence of Freedman's Classification of topological 4--manifolds we
have:
\begin{thm}\label{t:homeo}
The smooth 4--manifolds $X_1,X_2, X_3$ are homeomorphic to 
the rational surface $\cphat$. 
\end{thm}
\begin{proof}
Since the complements of the configurations are simply connected and
the fundamental group $\pi _1(\partial B_{p,q})$ surjects onto the
fundamental group of $B_{p,q}$, simple connectivity of $X_1, X_2,
X_3$ follows from Van Kampen's theorem.  Computing the Euler
characteristics and signatures of these 4--manifolds, Freedman's Theorem
\cite{Fr} implies the statement.
\end{proof}

\subsection{A further example}
A slightly different construction can be carried out as follows.

\begin{lem}
The plumbing 4--manifold $C_{32,15}$ embeds into $\cpk \# 16 \cpkk$.
\end{lem}
Recall that $C_{32,15}$ is equal to the 4--manifold defined by the linear
plumbing with weights $(-2,-9,-5,-2,-2,-2,-2,-2,-2,-3)$.
\begin{proof}
We start again with a fibration $\eegy \to {\mathbb {CP}}^1$ with a
singular fiber of type $\III^*$, three fishtails and two sections, as
shown by Figure~\ref{f:singfib}. After blowing up the double points of
the three fishtail fibers, blow up at $A_1,A_2$, smooth the
intersections at $B_1,B_2,A_3$ and keep the transverse intersection at
$B_3$.  One further blow-up as it is described by
Figure~\ref{f:inside} (performed inside the dashed circle of
Figure~\ref{f:singfib}) and finally the blow-up of the transverse
intersection of the section $s_1$ with the singular fiber of type
$\III^*$ provides the desired configuration $C_{32,15}$ in $\cpk \# 16
\cpkk$.
\end{proof}
\begin{rem}
Consider the configuration of curves in $\cpk \#14
\cpkk$ given above without the two last blow-ups.  This configuration
provides a ``necklace'' of spheres as shown in Figure~\ref{f:neck}.
\begin{figure}[htb]
\begin{center}
\setlength{\unitlength}{1mm}
\includegraphics[height=3.5cm]{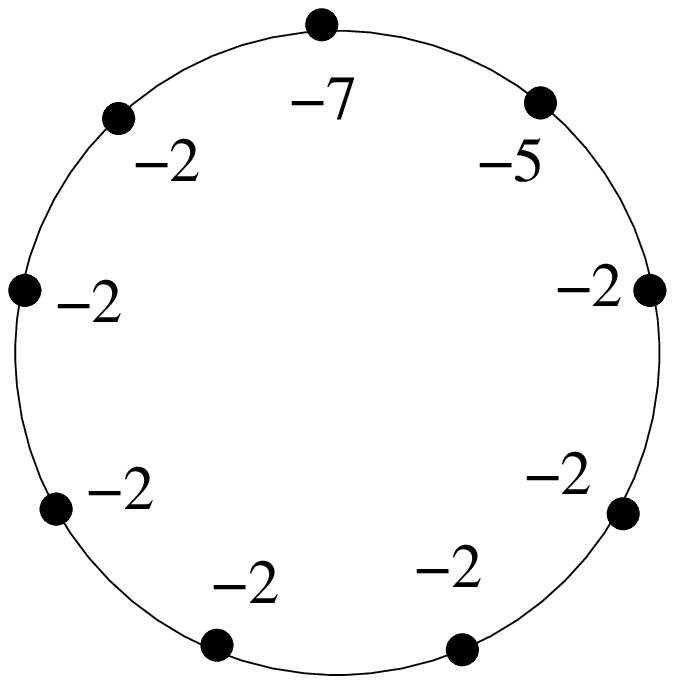}
\end{center}
\caption{``Necklace'' of spheres in $\cpk \# 14 \cpkk$}
\label{f:neck}
\end{figure}
Now $C_{32,15}$ can be given from this picture by blowing up the
intersection of the $(-7)$-- and the $(-2)$--framed circles, and then
blowing up the resulting $(-8)$--sphere appropriately one more time.
Notice that by blowing up the intersection of the $(-7)$-- and the
$(-5)$--curves instead, we can get two disjoint configurations of
$(-8, -2,-2, -2, -2)$ and $(-6,-2,-2)$, ie, two ``classical rational
blow-down'' configurations. Blowing them down we would recover the
existence of an exotic smooth structure on $\cphet$.
\end{rem}
Define $Y$ as the generalized rational blow-down of $\cpk \# 16 \cpkk$
along the configuration $C_{32,15}$ specified above.
\begin{lem}
The  4--manifold $Y$ is homeomorphic to $\cphat$.
\end{lem}
\begin{proof} 
Let $\alpha$ be the homology element represented by the circle in
$\partial C_{32,15}$ we get by intersecting the boundary of the
neighborhood of the plumbing with the sphere $S$ in the $\III^*$ fiber
not used in the construction, cf Figure~\ref{f:singfib}. Clearly
$\alpha$ is 9 times the normal circle of the last sphere in the
configuration. It follows that this circle generates $\pi_1 (\partial
C_{32,15})$.  Since it also bounds a disk in $\cpk \# 16 \cpkk -
{\mbox { int }}C_{32,15}$, the complement of the configuration is
simply connected, and so Van Kampen's theorem and the fact that 
the fundamental group $\pi _1(\partial B_{32,15})$ surjects onto the
fundamental group of $B_{32,15}$ shows simple connectivity. As before,
the computation of the Euler characteristics and signature, together
with Freedman's Theorem provides the result.
\end{proof}

\section{Seiberg--Witten invariants}
\label{s:sw}

In order to prove Theorem~\ref{t:main}, we will compute the
Seiberg--Witten invariants of the 4--manifolds constructed above.  In
order to make our presentation complete, we briefly recall basics of
Seiberg--Witten theory for 4--manifolds with $b_2^+=1$.  (For a more
thorough introduction to Seiberg--Witten theory, with a special
emphasis on the case of $b_2^+=1$, see \cite{Pratb, P, taubi}.) 

Suppose that $X$ is a simply connected, closed, oriented 4--manifold
with $b_2^+>0$ and odd, and fix a Riemannian metric $g$ on $X$.  Let
$L\to X$ be a given complex line bundle with $c_1(L)\in H^2(X; \bfz )$
characteristic, ie, $c_1(L)\equiv w_2(TX)$ (mod 2). Through its
first Chern class, the bundle $L$ determines a spin$^c$ structure $\s$
on $X$. The associated spinor $U(2)$--bundles $W^{\pm }_{\s}$ satisfy
$L\cong \det (W^{\pm}_{\s})$. A connection $A\in {\mathcal {A}}_L$ on
$L$, together with the Levi--Civita connection on $TX$ and the
Clifford multiplication on the spinor bundles induces a twisted Dirac
operator
\[
D_A\colon \Gamma (W^+_{\s})\to \Gamma (W^-_{\s}).
\]
For a connection $A\in {\mathcal  {A}}_L$, section $\Psi \in
\Gamma  (W^+_{\s})$ and  $g$--self--dual 2--form  $\eta\in \Omega  ^+_g (X;
\bfr )$ consider the \emph{perturbed Seiberg--Witten equations}
\[
D_A\Psi = 0, \qquad F_A^+=i(\Psi \otimes \Psi ^*)_0+i\eta,
\]
where $F_A^+$ is the self--dual part of the curvature $F_A$ of the
connection $A$ and $(\Psi \otimes \Psi ^*)_0$ is the frace--free part
of the endomorphism $\Psi \otimes \Psi ^*$.  For generic choice of the
self--dual 2--form $\eta$ the Seiberg--Witten moduli space --- which
is the quotient of the solution space to the above equations under the
action of the gauge group ${\mathcal {G}}={\rm {Aut}}(L)={\rm
{Maps}}(X; \bfr)$ --- is a smooth, compact manifold of dimension
\[
d_L=\frac{1}{4}(c_1^2(L)-3\sigma (X)-2\chi (X))
\]
(provided $d_L\geq 0$).  By fixing a 'homology orientation' on $X$,
that is, orienting $H^2_+(X; \bfr )$, the moduli space can be equipped
with a natural orientation. A natural 2--cohomology class $\beta $ can
be defined in the cohomology ring of the moduli space, and by
integrating $\beta ^{\frac{d_L}{2}}$ on the fundamental cycle of the
moduli space we get the \emph{Seiberg--Witten invariant} $SW_{X,g,\eta
}(L)$. This value is independent of the choice of the metric $g$ and
perturbation 2--form $\eta$ provided the manifold $X$ satisfies
$b_2^+(X)>1$. In case of $b_2^+(X)=1$, however, this independence
fails to hold. Let $\omega _g$ denote the unique self--dual 2--form
inducing the chosen homology orientation. It can be shown that $SW_{X,
g,h}(L)$ depends only on the sign of the expression
\[
(2\pi c_1(L)+[\eta ])\cdot [\omega _g].
\]
By fixing a sign for the above expression we say that we fixed a
\emph{chamber} for $L$, and going from one chamber to the other we
cross a \emph{wall}. It has been shown \cite{L} that by crossing a
wall the value of the Seiberg--Witten invariant changes by $\pm 1$.  To
specify a chamber, we need to fix a cohomology class (or its
Poincar\'e dual) with nonnegative square which can play the role of
$[\omega _g]$ for some metric. To remove the ambiguity on sign, we
require this element to pair positively with the element representing
the given homology orientation.

It is not hard to see that if $b_2^+(X)=1$ and $b_2^-(X)\leq 9$ then
$d_L\geq 0$ implies that $c_1^2(L)\geq 0$, hence for choosing the
perturbation term $\eta $ small in norm, the sign of $(2\pi
c_1(L)+[\eta ])\cdot [\omega _g]$ will be independent of the choice of the
metric $g$. Consequently, by restricting ourselves to Seiberg--Witten
invariants with small perturbation, on a manifold $X$ homeomorphic to
$\cpk \# n \cpkk$ with $n \leq 9$ the function
\[
SW_X \colon H^2(X; \bfz )\to \bfz
\]
is a diffeomorphism invariant.  For such a manifold a cohomology class
$K\in H^2 (X; \bfz )$ is called a \emph{Seiberg--Witten basic class}
if $SW_X(K)\neq 0$.  

It is a standard fact that, because of the presence of a metric with
positive scalar curvature, the Seiberg--Witten map vanishes for the
smooth 4--manifolds $\cpk \# n \cpkk$ with $n\leq 9$. Therefore in
order to show that the manifolds $X_i$ ($i=1,2,3$) given in
Definition~\ref{d:ex} provide exotic structures on $\cpk \# 6 \cpkk$
we only need to show that $SW_{X_i}\neq 0$.  We will go through the
computation of the invariants of $X_1$ only, the other cases follow
similar patterns.
\begin{thm}\label{t:sw}
There is a characteristic cohomology class ${\tilde {K}}\in H^2(X_1;\Z )$ with 
$SW_X({\tilde {K}} )\neq 0$.
\end{thm}
\begin{cor}\label{t:nondiffeo}
The 4--manifold $X_1$ is not diffeomorphic to $\cphat$.
\end{cor}
\begin{proof}
The corollary easily follows from Theorem~\ref{t:sw}, together with
the facts that $SW_{\cphat}\equiv 0$ and that the Seiberg--Witten
function is a diffeomorphism invariant for manifolds homeomorphic to
$\cpk \# n \cpkk$ with $n \leq 9$.
\end{proof}
\begin{proof}[Proof of Theorem~\ref{t:sw}]
Let $K\in H^2 (\cpk \# 17\cpkk ; \Z )$ denote the characteristic
cohomology class which satisfies
\[
K(h)=3 \qquad {\mbox {and}} \qquad K(e_i)=1\ (i=1, \ldots , 17).
\]
(The Poincar\'e dual of $K$ is equal to $3h-\sum _{i=1}^{17}e_i$.)  It
can be shown that the restriction $K\vert _{\cpk\# 17 \cpkk - {\mbox{
int }}C_{28,9}}$ extends as a characteristic cohomology class to $X_1$:
if $\mu $ denotes the generator of $H_1(\partial C_{28,9}; \Z )$
which is the boundary of a normal disk to the left--most circle in the
plumbing diagram of Figure~\ref{f:plum} then
\[
PD(K\vert _{\partial C})=532 \mu = 19 \cdot (28\mu)\in H_1(\partial
C_{28,9};\Z );
\]
since $ H_1 (B_{28,9}; \Z )$ is of order 28, the extendability
trivially follows. (See also Proposition~\ref{p:ch}.)  Let ${\tilde
{K}}$ denote the extension of $K\vert _{\cpk \# 17\cpkk - {\mbox { int
}}C_{28,9}}$ to $X_1$.  Using the gluing formula for Seiberg--Witten
invariants along lens spaces, see eg \cite{FS1, Pratb}, and the fact
that the dimensions of the moduli spaces defined by $K$ and ${\tilde
{K}}$ are equal, we have that the invariant $SW_{X_1}({\tilde {K}})$
is equal to the Seiberg--Witten invariant of $\cpk \# 17 \cpkk$
evaluated on $K$, in the chamber corresponding to a metric which we
get by pulling out $C_{28,9}$ along the 'neck' $L(784,251)\times
[-T,T]$ far enough.  For such a metric the period point provided by
the harmonic 2--form $\omega _g$ will be orthogonal to the
configuration $C_{28,9}$, hence the chamber can be represented by the
Poincar\'e dual of any homology element $\alpha\in H_2(\cpk \# 17\cpkk
; \Z)$ of nonnegative square represented in ${\cpk\# 17 \cpkk -
{\mbox { int }}C_{28,9}}$. For example,
\[
\alpha =
7h-2e_1-3e_2-\sum _3 ^9 2e_i- e_{10}-e_{12}-2e_{13}-e_{16}-e_{17}
\]
is such an element. (Simple computation shows that $\alpha $ is
orthogonal to all second homology elements in $C_{28,9}$, $\alpha \cdot
\alpha =0$ and $ \alpha \cdot h=7$.)

It is known that in the chamber corresponding to $PD(h)$ the Seiberg--Witten
invariant of $\cpk \# 17 \cpkk$ vanishes, since this is the chamber
containing the period point of a positive scalar curvature metric,
which prevents the existence of Seiberg--Witten solutions. Since the
wall--crossing phenomenon is well--understood in Seiberg--Witten
theory (the invariant changes by one once a wall is crossed), the
proof of the theorem reduces to determine whether $PD(\alpha)$ and
$PD(h)$ are in the same chamber with respect to $K$ or not.  Since
$K(h)=3>0$ and $h\cdot \alpha >0$, the inequality $K(\alpha )<0$ would
imply the existence of a wall between $PD(h)$ and $PD(\alpha) $, hence
$SW_{X_1} ({\tilde {K}})=SW_{\cpk \# 17 \cpkk }(K)\neq 0$, where the
invariant of $\cpk \# 17\cpkk$ is computed in the chamber containing
$PD(\alpha)$.  Simple computation shows that $K(\alpha )=-4$,
concluding the proof.
\end{proof}

\begin{proof}[Proof of Theorem~\ref{t:main}]
Now Theorem~\ref{t:homeo} and Corollary~\ref{t:nondiffeo} provide a
proof of the main theorem of the paper.
\end{proof}

In fact, with a little more effort we can determine all the Seiberg--Witten 
basic classes of $X_1$: 

\begin{prop}\label{p:min}
If $L\in H^2(X_1; \Z )$ is a Seiberg--Witten basic class of $X_1$ then
$L$ is equal to $\pm {\tilde {K}}$. Consequently $X_1$ is a minimal
4--manifold. 
\end{prop}

\begin{proof}
We start by studying $H_2(X_1-B_{28,9}; \Z)=
H_2(\cpk \# 17 \cpkk - C_{28,9}; \Z)$.
Clearly this is given by the subgroup of elements of 
$ H_2(\cpk \# 17 \cpkk ; \Z)$ that have trivial intersection in homology
with all the spheres in $C_{28,9}$.
A quick computation gives the following basis for this subgroup:
\begin{gather*}
A_1=e_{13}-e_{14},\ A_2=e_{14}-e_{15},\ A_3=e_{11}-e_{12},\ A_4=3h-e_{13}-
\sum _1 ^9 e_i,\\ 
A_5=-2e_1+2e_2-e_{11},
\ A_6=4h-e_1-2e_2-\sum _3 ^9 e_i-2e_{11}-e_{12}-e_{13},\\ 
A_7=h-e_2-e_{10}-e_{11}-e_{16}-e_{17}.
\end{gather*} 
Let $L$ be a Seiberg--Witten basic
class of $X_1$; then $L$ is uniquely determined by its restriction $L'
\in H^2(X_1-B_{28,9};\Z )$.  Following the argument in \cite{OSZP} we
determine the basic classes of $X_1$ in two steps.

First we select some smoothly embedded spheres and tori
in $\cpk \# 17 \cpkk$ that have trivial intersection number in homology
with all the spheres in the embedded configuration $C_{28,9}$. To this end note
that $A_1$, $A_2$, $A_3$, $A_5$,  $A_7$  can be represented by spheres and 
$A_4$, $A_6$ by tori. In addition, we will also use the
classes $A_8=A_1+A_4$ and $A_9=A_1+A_2+A_4$ --- these classes
 can be represented by tori.
In this first round we only search for basic classes 
$L$ that satisfy the additional
 adjunction inequalities
\begin{equation} \label{e:adj}
(A_i)^2+ |L(A_i)|\leq 0
\end{equation}
for $1\leq i\leq 9$. $L$ is determined by its evaluation on
  $A_1,\ldots ,A_7$, so the adjunction inequality on these elements
  leaves 8100 characteristic classes $L' \in H^2(X_1- B_{28,9};\Z)$ to
  consider. By the dimension formula, for a Seiberg--Witten basic
  class of $X_1$ we have $L^2=(L')^2 \geq 3$, and $(L')^2\equiv 3\
  {\rm mod}\ 8$. This test weeds out most of the classes: among the
  8100 classes there are only 22 with the right square.  Among these
  22 there are 20 that violate the adjunction inequality~(\ref{e:adj})
  along $A_8$ or $A_9$.  The remaining 2 classes evaluate as $\pm
  (0,0,0,1,1,2,2)$ on $A_1,\ldots ,A_7$ and thus correspond to $\mp
  K$.

To finish the computation let us assume that there is a
Seiberg--Witten basic class $L$ of $X_1$ that violates the adjunction
inequality~(\ref{e:adj}) with one of $A_i$. Note that any spin$^c$
structure on $\partial C_{28,9}$ that extends to $B_{28,9}$ has an
extension to $C_{28,9}$ with square equal to
$-b_2(C_{28,9})=-11$. Using such an extension and the gluing formula
along the lens space $\partial C_{28,9}$ we get a Seiberg--Witten
basic class $L_1$ of $\cpk \# 17 \cpkk$ in a chamber perpendicular to
the $C_{28,9}$ configuration, satisfying $d_L=d_{L_1}$ (where $d_L,
d_{L_1}$ denote the formal dimensions of the Seiberg--Witten moduli
spaces). Now using the adjunction relation for spheres and tori of
negative self--intersection \cite{FS2, OSZthom} we get another basic
class $L_2$ of $\cpk \# 17 \cpkk$ in a similar chamber with
$d(L_2)>d(L_1)$. By the gluing formula again, $L_2$ gives rise to a
basic class $L_3$ of $X_1$ with $d(L_3)>d(L)$.  Since we consider
Seiberg--Witten invariants with small perturbation term only, $X_1$
has a unique chamber. Therefore it has only finitely many basic
classes, consequently the above process has to stop, see also
\cite{OSZP}. It can stop only at a basic class that satisfies all the
adjunction inequalities for embedded spheres and tori and has positive
formal dimension.  Since $d_K=d_{-K}=0$, our previous search rules
this case out.
\end{proof}

\begin{rem} 
A similar computation applied to $X_2, X_3$ and $Y$ provides\break the same
result, hence these manifolds are also minimal, homeomorphic to\break
$\cphat$ but not diffeomorphic to it.  In particular, Seiberg--Witten
invariants do not distinguish these 4--manifolds from each other. 
\end{rem}

\begin{proof}[Proof of Corollary~\ref{hanykul}]
According to Proposition~\ref{p:min} and  \cite[Theorems~1.1, 1.2]{OSZP},
there are  4--manifolds $X,P,Q$
homeomorphic to $\cpk \# n \cpkk$ with $n=6,7,8$, resp.,
admitting exactly two Seiberg--Witten basic classes. Blowing up $X$ in at
most two and $P$ in at most one point, the application of the
blow-up formula for Seiberg--Witten invariants implies the corollary.
\end{proof}

\section{Symplectic structures}

Since our operation is a special case of the generalized rational blow-down
process, which is proved to be symplectic when performed along 
symplectically embedded spheres \cite{Sym}, we conclude:

\begin{thm}
The 4--manifolds $X_1, X_2, X_3$ and $Y$ constructed above admit
symplectic structures.
\end{thm}
\begin{proof}
The 2--spheres in the configurations are either complex submanifolds or
given by smoothings of transverse intersections of complex submanifolds,
which are known to be symplectic. Furthermore, all geometric intersections
are positive, hence the result of \cite{Sym} applies.
\end{proof}

In the rest of the section we study the limit of Seiberg--Witten
invariants in detecting exotic smooth 4--manifolds with symplectic structures which are 
homeomorphic to small rational surfaces.

\begin{prop}\label{p:mini}
Suppose that the smooth 4--manifold $X$ is homeomorphic to $S^2\times
S^2$ or $\cpk \# n \cpkk$ with $n \leq 8$ and $X$ admits a symplectic
structure $\omega$. If $X$ has more than one pair of Seiberg--Witten
basic classes then $X$ is not minimal.
\end{prop}
\begin{proof}
By \cite{L} we know that if $c_1(X)\cdot [\omega ] >0$ and $X$ is
simply connected then $X$ is a rational surface, hence under the above
topological constraint it admits no Seiberg--Witten basic classes.
Therefore we can assume that $c_1(X)\cdot [\omega ]<0$.  Suppose now
that $\pm K$ and $\pm L$ are both pairs of basic classes and $K\neq
\pm L$.  Notice that by the dimension formula for the Seiberg--Witten
moduli spaces it follows that $K^2>0$ and $L^2>0$.
Suppose furthermore that $X$ is minimal.
By a theorem of Taubes \cite{Taub} we can assume that
$K=-c_1(X)$ and we can choose the sign of $L$ to satisfy $L\cdot
[\omega ]>0$.  Let $a$ denote the Poincar\'e dual of the cohomology
class $\frac{1}{2}(K-L)$.  By \cite{Taubmas} and the fact that
$SW_X(-L)\neq 0$, the nontrivial homology class $a$ can be represented
by a pseudo--holomorphic curves. It follows then that $(K-L)\cdot
[\omega ]>0$.

Suppose first that $(K-L)^2\geq 0$. 
Then the Light Cone Lemma \cite[Lemma~2.6]{L} implies that
$K\cdot (K-L)>0$ and $L\cdot (K-L)>0$ unless $L=rK$ for some $r\in \bfq$.
The two inequalities imply $K^2>L^2$, contradicting the fact that the
moduli space corresponding to the spin$^c$ structure determined by $L$ is
of nonnegative formal dimension. If $L=rK$ and $K\cdot (K-L)=0$, then 
the fact $K^2>0$ implies that $r=1$, hence $L=K$, contradicting
our assumption $K\neq \pm L$.

Finally we have to examine the case when $(K-L)^2<0$.  In this case,
by \cite[Proposition~7.1]{Taubmas}
for generic almost--complex structure the pseudo--holomorphic
representative of the homology class $a = PD(\frac{1}{2}(K-L))$
contains a sphere component of square $(-1)$,
contradicting the minimality of $X$.
\end{proof}

\begin{prop}\label{mini2}
Suppose that $(X_1, \omega _1)$ and 
$(X_2, \omega _2)$ are simply connected minimal symplectic
 4--manifolds with $b_2^+=1$ and $b_2^- \leq 8$. If $X_1$ and $X_2$
 are homeomorphic and have nonvanishing Seiberg--Witten invariants
 then we can choose the homeomorphism $f\colon X_1\to X_2$ so that
$$ SW_{X_2}(L)=\pm SW_{X_1}(f^{\ast}(L))$$
for all characteristic classes $L\in H^2(X_2;\Z)$. 
\end{prop}

\begin{proof}
According to Proposition \ref{p:mini} both $X_1$ and $X_2$ has two
basic classes $\pm c_1(X_i, \omega _i)$. According to Taubes' theorem
\cite{Taub} (using an appropriate homology orientation) we have
\[
SW_{X_i}(c_1(X_i, \omega _i))=1, \ \ SW_{X_i}(-c_1(X_i, \omega _i))=-1.
\]
According to Freedman \cite{Fr} the required homomorphism $f$ can be
induced by an isomorphism
\[
g\co  H^2(X_2;\Z)\longrightarrow H^2(X_1;\Z)
\] 
that maps $c_1(X_2,\omega _2)$ to $c_1(X_1, \omega _1)$ and preserves
the intersection form.  The existence of such $g$ is trivial when
$b_2^-$ is zero or one; the general case follows from the large
automorphism group of the second cohomology group $H^2$ given by
reflecting on cohomology classes with squares $1$, $-1$ and $-2$. In
particular, for the intersection form of $\cpk \# n\cpkk$ ($2\leq
n\leq 8$) it is easy to use reflections along the Poincar\'e duals of
$h$, $e_i$, $h-e_i-e_j$ and $h-e_i-e_j-e_k$ to map a given
characteristic class $L$ with $L^2=9-n$ to $3h-e_1- \ldots -e_n$.
Depending on whether $g$ respects the chosen homology orientations on
$X_1$ and $X_2$ or not, we have $SW_{X_2}(L)=\pm
SW_{X_1}(f^{\ast}(L)).$
\end{proof}

The above result together with the blow-up formula for Seiberg--Witten
invariants imply the following:

\begin{cor}
The Seiberg--Witten invariants can distinguish at most finitely many
symplectic 4--manifolds homeomorphic to a rational surface $X$
with Euler characteristic $e(X)<12$. \qed
\end{cor}

\section{Appendix: singular fibers in elliptic fibrations}
\label{app}

For the sake of completeness we give an explicit construction of the
elliptic fibration $\eegy \to {\mathbb {CP}}^1$ used in the paper.
The existence of such fibration is a standard result in complex
geometry; in the following we will present it in a way useful for
differential topological considerations.

Notice first that to verify the existence of a fibration with singular fibers
of type $\III^*$ (also known as the ${\tilde {E}}_7$--fiber) and three fishtail 
fibers is quite easy. As it is shown in \cite{Kod} (see also 
\cite[pp. 35--36]{HKK})
the monodromy of an ${\tilde {E}}_7$--fiber can be chosen to 
be equal to  $\left(\begin{smallmatrix}
 0 & -1 \\
 1 & 0 
\end{smallmatrix}
\right)$, while for a fishtail fiber the monodromy is conjugate to 
$\left(\begin{smallmatrix}
 1 & 1 \\
 0 & 1
\end{smallmatrix}
\right)$.  
Since
\[
\left(\begin{matrix}
 0 & -1 \\
 1 & 0 
\end{matrix}
\right) 
\left(\begin{matrix}
 1 & 1 \\
 0 & 1
\end{matrix}
\right)
\left(\left(\begin{matrix}
 0 & 1 \\
 -1 & 0
\end{matrix}
\right)^{-1}
\left(\begin{matrix}
 1 & 1 \\
 0 & 1
\end{matrix}
\right)
\left(\begin{matrix}
 0 & 1 \\
 -1 & 0
\end{matrix}
\right)\right)
\left(\begin{matrix}
 1 & 1 \\
 0 & 1
\end{matrix}
\right)=
\left(\begin{matrix}
 1 & 0 \\
 0 & 1
\end{matrix}
\right) ,
\]
the genus--1 Lefschetz fibration with the prescribed singular fibers
over the disk extends to a fibration over the sphere $S^2$. Simple
Euler characteristic computation and the classification of genus--1
Lefschetz fibrations show that the result is an elliptic fibration on
$\eegy$. The existence of the two sections positioned as required by
the configuration of Figure~\ref{f:singfib} is, however, less apparent
from this picture. One could repeat the above computation in the
mapping class group of the typical fiber with appropriate marked
points, arriving to the same conclusion.  Here we rather use a more
direct way of describing a pencil of curves in $\cpk$ and following
the blow-up procedure explicitly.

\begin{figure}[htb]
\begin{center}
\setlength{\unitlength}{1mm}
\includegraphics[height=3.5cm]{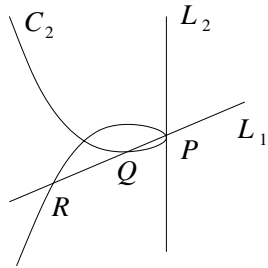}
\end{center}
\caption{Curves generating the pencil}
\label{f:curves}
\end{figure}

Let 
\[ 
C_1=\{ [x:y:z]\in \cpk \mid p_1(x,y,z)=(x-z)z^2=0\} \qquad {\mbox { and }}
\]
\[
C_2=\{ [x:y:z]\in \cpk \mid p_2(x,y,z)=x^3+zx^2-zy^2=0\} 
\]
be two given complex curves in the complex projective plane $\cpk$.
The curve $C_1$ is the union of the lines $L_1=\{ (x-z)=0\}$ and
$L_2=\{ z=0\}$, with the latter of multiplicity two. $C_2$ is an
immersed sphere with one positive transverse double point --- blowing
this curve up nine times in its smooth points results a fishtail
fiber, see also \cite[Section~2.3]{GS}. $L_2$ intersects $C_2$ in a
single point $P=[0:1:0]$ (hence this point is a triple tangency
between the two curves), and $L_1$ (also passing through $P$)
intersects $C_2$ in two further (smooth) points $R=[1:\sqrt{2}:1]$ and
$Q=[1:-\sqrt{2}:1]$, cf Figure~\ref{f:curves}. Therefore the
pencil 
\[ 
C_t=C_{[t_1:t_2]}=\{ (t_1p_1+t_2p_2)^{-1}(0)\} \qquad (t=[t_1:t_2]\in
{\mathbb {CP}}^1)
\]
of elliptic curves defined by $C_1$ and $C_2$
provides a map $f$ from $\cpk $ to ${\mathbb {CP}}^1$ well--defined
away from the three base points $P,Q,R$. In order to get the desired
fibration we will perform seven infinitely close blow-ups at the
base point $P$ and two further blow-ups at $R$ and $Q$, resp.
We will explain only the first blow-up at $P$, the rest 
follows a similar pattern.
After the blow-up of $P$ we would like to have a pencil on the
blown-up manifold. We take ${\tilde {C}}_2$ to be the proper
transform of $C_2$, while ${\tilde {C}}_1$ will be the proper
transform of $C_1$ together with a certain multiple of the exceptional
divisor, chosen so that the two curves represent the same homology
class. Under this homological condition the two curves can be given as
zero sets of holomorphic sections of the same holomorphic line bundle,
hence ${\tilde {C}}_1$ and ${\tilde {C}}_2$ define a pencil on the
blown-up manifold.  Since $[C_1]=[C_2]=3h \in H_2 (\cpk ; \Z)$, it is
easy to see that $[{\tilde {C}}_2]=3h-e_1 \in H_2(\cpk \# \cpkk ;
\Z)$, where (as usual) $e_1$ denotes the homology class of the
exceptional divisor of the blow-up. Now it is a simple matter to see
that the proper transform of $C_1$ is the union of the proper
transforms of $L_1$ and $L_2$, which transforms represent $h-e_1$ and
$2h-2e_1$ in $H_2(\cpk \# \cpkk ; \Z)$. Therefore in the pencil we
need to take the curve given by the proper transform of $C_1$ together
with the exceptional curve, the latter with multiplicity two.  Since
$e_1$ is part of ${\tilde {C}}_1$, we further have to blow up its
intersection with ${\tilde {C}}_2$.  The same principle shows that in
the further blow-ups the exceptional divisors $e_2,e_3,e_4, e_5, e_6$
come with multiplicities $3,4,3,2$ and 1. Finally, after blowing up
for the seventh time, the two curves defining the pencil get locally
separated, and hence $e_7$ will not lie in any of the curves of the
new pencil anymore --- it will be a section, that is, it intersects
all the curves in the pencil transversally in one point. (Notice that
we used seven blow-ups to separate the curves $C_1$ and $C_2$ at $P$,
where they intersected each other of order seven: $L_2$ being a linear
curve of multiplicity two, intersected the cubic curve $C_2$ of order
six, while $L_1$ simply passed through the intersection point $P$.)
After blowing up the two further base points $Q,R$, we get a fibration
on the nine--fold blow-up $\eegy$ with a fishtail fiber, a singular
fiber of type $\III^*$ and three sections provided by the exceptional
divisors $e_7,e_8$ and $e_9$.  (To recover the homology classes of the
spheres in the type $\III^*$ fiber indicated by Figure~\ref{f:III} one
needs to rename the exceptional divisors of the blow-ups; we leave
this simple exercise to the reader.) A final simple calculation shows
that the resulting fibration has three fishtail fibers:

\begin{prop} The pencil 
\[ 
\{ C_{[t_1:t_2]}= (t_1 p_1+t_2p_2)^{-1}(0)\mid [t_1:t_2]\in  
{\mathbb {CP}}^1\} 
\] 
contains four singular curves: $C_1, C_2$ and $C_3, C_4$. Furthermore
the latter two curves are homeomorphic to $C_2$ and give rise to 
fishtail fibers after blowing up the base points of the pencil.
\end{prop}
\begin{proof}
Since $L_2=\{ z=0\} \subset C_1$, all other curves of the pencil
are contained in $\{ z=1\}\cup \{ P\}$. The curve 
$C_t=C_{[t_1:t_2]}$ has a singular point if and 
only if for the polynomial
$p_t(x,y)=t_1(x-1)+t_2(x^3+x^2-y^2)$
we can find  $(x_0,y_0)\in \C ^2$ with
\[
p_t(x_0,y_0)=0, \qquad \frac{\partial p_t}{\partial x}(x_0,y_0)=0 \qquad
{\mbox {  and  }} \qquad\frac{\partial p_t}{\partial y}(x_0,y_0)=0.
\]
Since $ \frac{\partial p_t}{\partial y}(x,y)=-2t_2y$, it vanishes
if $t_2=0$ (providing $C_1$) or $y=0$.
In the latter case the above system reduces to
$$
t_1(x-1)+t_2(x^3+x^2)=0 \qquad {\mbox{ and }} \qquad t_1+t_2(3x^2+2x)=0,
$$
which admits a nontrivial solution $(t_1, t_2)$ if and only if the
determinant 
\[ 
(x-1)(3x^2+2x)-(x^3+x^2)=2x(x^2-x-1)
\]
vanishes. The solution $x=0$ implies $t_1=0$, giving $C_2$, while 
$x=\frac{1}{2}(1\pm \sqrt{5})$ give the two singular points on the curves 
$C_3$ and $C_4$. Simple Euler characteristic computation
shows that the two curves will give rise to fishtail fibers in the 
elliptic fibration.
\end{proof}

\end{document}